\documentclass[12pt,oneside]{amsart} 
\usepackage{amssymb,amscd}
\usepackage{euscript}
 
      \theoremstyle{plain}
      \newtheorem{theorem}{Theorem}[section]
      \newtheorem{lemma}[theorem]{Lemma}
      \newtheorem{corollary}[theorem]{Corollary}
      \newtheorem{proposition}[theorem]{Proposition}

      \newtheorem{definition}[theorem]{Definition}        
    \newtheorem{assumptions}[theorem]{Standing Assumptions}     
          
\numberwithin{equation}{section}

      \makeatletter
      \def\@setcopyright{}
      \def\serieslogo@{}
      \makeatother

\def\A{\EuScript{A}} 
\def\B{\EuScript{B}} 
 
\def\E{\mathcal{V}}
\def\V{\mathcal{V}}

\def\M{X}
\def\m{\mathcal{M}}

\def\R{\mathbb R}
\def\Z{\mathbb Z}
\def\N{\mathbb N}

\def\dist{\text{dist}}

\def\Id{\text{Id}}
\def\e{\epsilon}

\def\la{\lambda}

\def\bv{\mathbf v}

\def\QED{\hfill\hfill{\square}}

\begin{document}

\author{Victoria Sadovskaya$^{\ast}$}

\address{Department of Mathematics, The Pennsylvania State University, University Park, PA 16802, USA.}
\email{sadovskaya@psu.edu}

\title[Fiber bunching and cohomology  for  Banach cocycles]
{Fiber bunching and cohomology  for \\ Banach cocycles over hyperbolic systems} 

\thanks{{\it Mathematical subject classification:}\,  37D20, 37C15, 37H05}
\thanks{{\it Keywords:} Cocycle, cohomology, hyperbolic system, periodic point, fiber bunching,
Banach space}
\thanks{$^{\ast}$ Supported in part by NSF grant DMS-1301693}


\begin{abstract} 

We consider H\"older continuous cocycles over hyperbolic dynamical systems
with values in the group of invertible bounded linear operators on a Banach space. 
We show that two fiber bunched cocycles are H\"older continuously cohomologous 
if and only if they have H\"older conjugate periodic data.  The fiber bunching condition 
means that non-conformality of the cocycle is dominated by the expansion and 
contraction in the base system. We show that this condition can be established 
based on the periodic data of a cocycle.  We also establish H\"older continuity of 
a measurable conjugacy between a fiber bunched cocycle and one with values 
in a set which is compact in strong operator topology.

\end{abstract}

\maketitle 

\section{Introduction and statements of the results}
Cocycles play an important role in dynamics. Cohomology of real-valued 
and, more generally, group-valued cocycles over hyperbolic systems has 
been extensively studied starting with the seminal work of A. Liv\v{s}ic \cite{Liv1},
see \cite{KtN} for an overview. The study has been focused on obtaining cohomology 
of two cocycles from their periodic data, i.e. the values at the periodic points of the 
base system, and on regularity of transfer map, or conjugacy, between two cocycles. 
 Liv\v{s}ic resolved the case of cocycles with values in $\R$ or an abelian group and 
made some progress for more general groups.
For smooth dynamical systems, the differential and its restrictions to invariant 
sub-bundles give important examples of cocycles. This motivated in part
the extensive research of $GL(n,\R)$ and Lie group valued cocycles.
Cohomology problems for cocycles with values in non-abelian groups are much more difficult.
The case when one of the cocycles is the identity has been studied most and by now is 
relatively well understood, see for example \cite{Liv2, NT95, PW, LW10, K11, GG}.
The cohomology problem for two arbitrary cocycles does not reduce to this special case 
for non-abelian groups. This problem was first considered in \cite{Pa} for compact groups and in
 \cite{Sch} for cocycles with ``bounded distortion". The most general results for $GL(n,\R)$ 
 were established in \cite{S15} for fiber bunched cocycles. One of these results was also
 independently obtained in \cite{B}.

The infinite dimensional case is more difficult and so far is much less developed. 
For two arbitrary cocycles with values in the group of diffeomorphisms of a compact 
manifold, higher regularity of the conjugacy was studied in \cite{NT98} and recently
cohomology of cocycles with equal periodic data was obtained in  \cite{BK} under
a certain  bunching assumption on both cocycles. 
In this paper we extend the results for finite dimensional linear cocycles to the 
infinite-dimensional setting.  
We consider cocycles with values in the group of invertible operators on a 
Banach space $V$.  The simplest examples  are given by random and Markov 
sequences of operators. They correspond to locally constant cocycles over
 subshifts of finite type.

The space  $L(V)$ of bounded linear operators on $V$ 
 is  a Banach space equipped 
with the operator norm $\|A\|=\sup \,\{ \|Av\| : \,v\in V , \;\|v\| \le 1\}.$ 
The open set $GL(V)$  of invertible elements in  $L(V)$ is a topological group  
and a complete metric space with respect  to the metric
\begin{equation}\label{d}
d (A, B) = \| A  - B \|  + \| A^{-1}  - B^{-1} \|.
\end{equation}

\begin{definition} Let $f$ be a homeomorphism of a compact metric space $X$
and let $A$ be a function from $X$ to $(GL(V),d)$. 
The {\em Banach  cocycle over $f$ generated by }$A$ 
is the map $\A:\,X \times \Z \,\to G$ defined  by $\,\A(x,0)=\Id\,$ and for $n\in \N,$
 $$
\A(x,n)=\A_x^n = A(f^{n-1} x)\circ \cdots \circ A(x) \quad\text{and}\quad\,
\A(x,-n)=\A_x^{-n}= (\A_{f^{-n} x}^n)^{-1} .
$$
Clearly, $\A$ satisfies the {\em cocycle equation}\,
$\A^{n+k}_x= \A^n_{f^k x} \circ \A^k_x$.
\vskip.05cm
We say that the cocycle $\A$ is   {\em $\beta$-H\"older} if its generator $A$  is H\"older continuous with exponent $0<\beta \le 1$ with respect to the metric $d$, i.e. there exists $K>0$ such that
$$
 d (A(x), A(y)) \le  \, K \dist (x,y)^\beta \quad\text{for all }x,y\in X.
$$
\end{definition}

H\"older continuity  is needed to develop a theory even for
scalar cocycles. On the other hand, higher regularity is rare
for cocycles given by restrictions 
of the differential of an Anosov diffeomorphism to the stable and unstable subbundles.
Additionally, symbolic dynamical systems
have a natural H\"older structure, but  not a smooth one.


\begin{definition}
The {\em quasiconformal distortion} of a cocycle $\A$ is the function
 $$
Q_\A(x,n)= \| \A_x^n\| \cdot \| (\A_x^n)^{-1}\|, 
 \quad x\in X \text{ and }n\in \Z.
 $$ 
 \end{definition}

The quasiconformal distortion is a measure of non-conformality of the cocycle.  
If $Q_\A(x,n)\le K$ for all $x$ and $n$ the cocycle is said to be uniformly 
quasiconformal, and if $Q_\A(x,n)=1$ for all $x$ and $n$ it is said 
to be conformal.

\begin{definition} \label{bunching def}
 A $\beta$-H\"older  cocycle $\A$ over a hyperbolic system $(X,f)$ as in Section~\ref{base} is  
 {\em fiber bunched} if 
there exist numbers $\theta<1$ and $L$  such that for all $x\in\M$ and $n\in \N$,
\begin{equation}\label{fiber bunched}
Q_\A(x,n) \cdot  (\nu^n_x)^\beta < L\, \theta^n \quad\text{and}\quad
Q_\A(x,-n) \cdot  (\hat \nu^{-n}_x)^\beta < L\, \theta^n, 
\end{equation}
where 
$ \;\nu^n_x=\nu(f^{n-1}x)\cdots\nu(x) \,\text{ and }\;
   \hat\nu^{-n}_x=(\hat \nu(f^{-n}x))^{-1}\cdots  (\hat\nu(f^{-1}x))^{-1}.$
\end{definition}

This condition means that non-conformality of the cocycle is dominated by the 
expansion and contraction in the base. In particular, bounded, conformal, and uniformly 
quasiconformal cocycles are fiber bunched. 
The condition guarantees convergence of certain iterates of the cocycle 
along the stable and unstable leaves, and it plays an important role in the study of 
non-commutative cocycles.
We use the weakest, ``pointwise", version of  fiber bunching with non-constant
estimates of expansion and contraction in the base.

\begin{assumptions} In the results below, $(X,f)$ is a hyperbolic system as described in 
Section \ref{base},
 and $\A$ and $\B$ are $\beta$-H\"older  Banach cocycles over $f$.
\end{assumptions}

The next proposition allows us to obtain fiber bunching of a cocycle $\A$ from fiber 
bunching of its periodic data.

\begin{proposition} \label{bunching per}
Suppose that for a cocycle $\A$ there exist numbers $\,\tilde\theta<1$ and  $\tilde L$ such that whenever 
$f^kp=p$, $\,k\in \N$,\, we have,
\begin{equation}\label{fiber bunched per}
Q_\A(p,k) \cdot  (\nu^k_p)^\beta < \tilde L\, \tilde\theta^k \quad\text{and}\quad
Q_\A(p,-k) \cdot  (\hat \nu^{-k}_p)^\beta < \tilde L\, \tilde \theta^k.
\end{equation}
Then $\A$ is fiber bunched.
\end{proposition}

It follows that a cocycle with periodic data conjugate to that of a 
fiber bunched cocycle via a bounded conjugacy is also fiber bunched.

\begin{definition} We say that Banach cocycles $\A$ and $\B$ have 
{\em conjugate periodic data} if  for  every periodic point $p=f^k p$
 there exists  $C(p)\in GL(V)$  such that 
 $$
 \B_p^k=C(p) \circ \A_p^k \circ C(p)^{-1}.
$$
\end{definition}

\begin{corollary} \label{bunching conj}
Suppose that  $\A$ is fiber bunched and $\B$ has conjugate periodic data
with bounded conjugacy, i.e. $ \max\,\{\|C(p)\|,\, \|C(p)^{-1} \| \} \le M,$
 where $M$ is a constant independent of $p$.
Then $\B$ is also fiber bunched.
\end{corollary}

A natural equivalence for cocycles is cohomology, i.e. existence of a conjugacy, 
which can be considered in various regularity classes. 

\begin{definition} \label{cohdef}
Cocycles $\A$ and $\B$ are 
(measurably, continuously)  {\em cohomologous} 
if there exists  a (measurable, continuous) 
function  $C:\M\to GL(V)$, called a {\em conjugacy} or
a transfer map between $\A$ and $\B$, such that
\begin{equation}\label{conj eq}
  \A_x^n=C(f^n x) \circ \B_x^n \circ  C(x)^{-1} 
  \quad\text{ for all }x\in \M  \text{ and }n\in \Z.
\end{equation}
\end{definition}

Clearly, if two cocycles are continuously cohomologous, then they have conjugate 
periodic data. The converse is not true in general\, even when $V$ is two-dimensional
and $C(p)$ is bounded  \cite{S13}. If $C(p)$ is  H\"older, conjugating 
$\B$ by the extension of $C$ reduces the problem to the case
of equal periodic data, i.e. $\A_p^n= \B_p^n$. Positive results for  
equal periodic data, as well as some results for conjugate data, 
were established by W.~Parry \cite{Pa} for compact $G$ and, 
somewhat more generally, by K. Schmidt \cite{Sch} for cocycles 
with ``bounded distortion". First results outside this setting were 
obtained in  \cite{S13} for certain types of $GL(2, \R)$-valued cocycles.
In \cite{S15} we considered $GL(n,\R)$-valued cocycles over hyperbolic systems.
We showed that if a cocycle $\A$ is fiber bunched and 
$\B$ has equal periodic data, then $\A$ and $\B$ are 
H\"older continuously cohomologous.
Moreover, we obtained H\"older cohomology under the assumption that 
$\A$ is fiber bunched,  $\B$ has conjugate periodic data
and the conjugacy $C(p)$ is $\beta$-H\"older continuous at a fixed point of $f$. 
In the following theorem we extend the results to  Banach cocycles and 
remove the assumption that $f$ has a fixed point. 

\begin{theorem} \label{conjugate data} 
Suppose that  a cocycle $\A$ is fiber bunched,  $\B$ has conjugate periodic data, 
and $C(p)$ is $\beta$-H\"older continuous  at a periodic point $p_0$, i.e.
there exists a constant $c$ such that 
$\,d(C(p), C(p_0))\le c\,\dist (p,p_0)^\beta$ for every periodic point $p$.

Then there exists a unique $\beta$-H\"older continuous conjugacy $\bar C$ 
between $\A$ and $\B$ such that $\bar C (p_0)=C(p_0)$. 
\end{theorem}

\noindent We do not assume fiber bunching for $\B$ as it  follows immediately from
Corollary \ref{bunching conj}.
\vskip.1cm

We note that  $\bar C(p)$ does  not necessarily coincide with 
$C(p)$ for $p\ne p_0$. For example, let $\B\equiv \Id$ and
let $\A_x=\bar C(fx) \circ \bar C(x)^{-1}=\bar C(fx) \circ \B_x\circ  \bar C(x)^{-1}$, 
where $C$ is a H\"older 
continuous function with $\bar C (p_0)=\Id$.
Then $\A^n_p=\B^n_p=\Id\,$ whenever $p=f^np$, and so we can take
$C(p)=\Id$ for each $p$.

\begin{corollary} \label{equal data} 
Suppose that a cocycle $\A$ is fiber bunched and a cocycle $\B$ has equal
 periodic data, i.e. $\A_p^k=\B_p^k$ whenever $p=f^kp$.
Then the cocycles are $\beta$-H\"older continuously cohomologous.

In particular, if $\B$ is a cocycle with $\,\B_p^k=\Id\,$  whenever $p=f^kp$, then $\B$ is 
$\beta$-H\"older continuously cohomologous to the identity cocycle.
\end{corollary}

The second part of this corollary was recently obtained by G. Grabarnik and M. Guysinsky
for cocycles with values in Banach algebras \cite{GG}.
\vskip.2cm 

Next we consider the question whether a measurable conjugacy between two cocycles 
is continuous, i.e. coincides with a continuous conjugacy on a set of full measure.
Measurability is understood with respect to a suitable measure, for example the measure 
of maximal entropy or the invariant volume. This problem was also first considered in
the case when one of the cocycles is the identity. The first result beyond this case was obtained
 by K. Schmidt \cite{Sch} for two cocycles with ``bounded distortion", the prime examples
 being uniformly bounded $GL(n,\R)$-valued cocycles and ones with values in compact groups.
A counterexample by M. Pollicott and C. P. Walkden  \cite{PW}  showed that  additional 
assumptions are needed for a positive answer in more general context: they
constructed $GL(2,\R)$-valued cocycles which are measurably but not continuously 
cohomologous. Moreover both cocycles can be made arbitrarily close to the identity, 
and in particular fiber bunched. In \cite{S15} we showed  that if $\A$ is fiber bunched 
and $\B$ is uniformly quasiconformal, i.e. $Q_\B(x,n)\le K$ for all $x$ and $n$,  then 
any measurable conjugacy between $\A$ and $\B$ is $\beta$-H\"older continuous. 

The following theorem extends the finite-dimensional results to the Banach setting. 
One of the difficulties here is that the space $(GL(V),d)$ is not separable even if $V$ is.
To use the tools of the theory of measurable functions, such as Lusin's theorem, 
we work with the strong operator topology, i.e. the topology of pointwise convergence.
We assume that $\B$ takes values in a precompact set, which for a finite dimensional 
$V$ is equivalent to uniform boundedness of $\B$ in $(GL(V),d)$.

\begin{theorem} \label{measurable conjugacy}
Suppose that the Banach space $V$ is separable, a cocycle $\A$ is fiber bunched and 
a cocycle $\B$ takes values in a subset of $GL(V)$ that is precompact in the topology 
of pointwise convergence. Let  $\mu$ be an ergodic
invariant  measure with full support and local product structure.
Then any $\mu$-measurable conjugacy between $\A$ and $\B$ coincides with a $\beta$-H\"older continuous conjugacy on a set of full measure.

\end{theorem}

A measure has local product structure if it is locally equivalent to the product of its conditional measures on the local stable and unstable manifolds.
Examples of ergodic measures with full support and local product structure include 
the measure of maximal entropy,  more generally Gibbs (equilibrium) measures of H\"older continuous potentials, and the invariant volume  if it exists.

\vskip.1cm


\section{Hyperbilic systems in the base} \label{base}

\noindent{\bf Transitive Anosov diffeomorphisms.} 
A diffeomorphism  $f$ of a compact connected  manifold $\M$
 is called {\it Anosov}\, if there exist a splitting 
of the tangent bundle $T\M$ into a direct sum of two $Df$-invariant 
continuous subbundles $E^s$ and $E^u$,  a Riemannian 
metric on $\M$, and  continuous  
functions $\nu$ and $\hat\nu$  such that 
\begin{equation}\label{Anosov def}
\|Df_x(v^s)\| < \nu(x) < 1 < \hat\nu(x) <\|Df_x(v^u)\|
\end{equation}
for any $x \in \M$ and unit vectors  
$\,\bv^s\in E^s(x)$ and $\,\bv^u\in E^u(x)$.
The sub-bundles $E^s$ and $E^u$ are called stable and unstable. 
They are tangent to the stable and unstable foliations 
$W^s$ and $W^u$ respectively. 

Using \eqref{Anosov def} we choose a 
small number $\rho>0$ such 
 that for every $x \in \M$ we have
 $\| Df_y\| < \nu (x)$ for all $y$ in the ball  in $W^{s}(x)$ centered 
at $x$ of  radius $\rho$  in the intrinsic metric of $W^{s}(x)$
 and such that $W^s_{\text{loc}}(x)\cap W^s_{\text{loc}}(z)$ consists of a single point 
for any sufficiently close $x$ and $z$ in $X$. This property is called the local 
product structure of the foliations.
We refer to this ball as the {\it local stable manifold} of $x$ and denote it by $W^{s}_{\text{loc}}(x)$.
Local unstable manifolds 
are defined similarly. It follows that for all $n\in \N$,
\begin{equation}\label{dist}
\begin{aligned}
&\dist (f^nx, f^ny)< \nu^n_x \cdot  \dist(x,y) \quad \text{for all } x \in \M
\text{ and }y\in W^s_{\text{loc}}(x),\\
&\dist (f^{-n}x, f^{-n}y)< \hat \nu^{-n}_x \cdot  \dist(x,y) \quad \text{for all } x \in \M
\text{ and }y\in W^u_{\text{loc}}(x).
\end{aligned}
\end{equation}

A diffeomorphism $f$ is {\it (topologically) transitive} if there is a point $x$ in $\M$
with dense orbit. All known examples of Anosov diffeomorphisms have this property.
\vskip.2cm


\noindent{\bf Topologically mixing diffeomorphisms of locally maximal hyperbolic sets.}  
More generally, let $f$ be a diffeomorphism of a manifold $\m$.
A compact $f$-invariant  set $X \subseteq \m$ is
called {\em hyperbolic} if there  exist a continuous $Df$-invariant splitting 
$T_X \m = E^s\oplus E^u$, and a Riemannian metric and 
continuous functions $\nu$, $\hat \nu$ on an open set
$U \supseteq X$ such that \eqref{Anosov def} holds for all $x \in X$.
The set $X$ is called {\em locally maximal} if 
$X= \bigcap_{n\in \Z} f^{-n }(U)$ for some open set $U\supseteq X$. 
 

\vskip.2cm
\noindent{\bf Mixing subshifts of finite type.}
Let $M$ be $k \times k$ matrix with entries from $\{ 0,1 \} $ such that all 
entries of $M^N$ are positive for some $N$. Let
\vskip.15cm
$\hskip1cm X= \{ \,x=(x_n) _{n\in \Z}\, : \,\; 1\le x_n\le k \;\text{ and }\;
 M_{x_n,x_{n+1}}=1 \,\text{ for every } n\in \Z \,\}.$
\vskip.15cm
\noindent The shift map $f:X\to X\,$ is defined by 
$(fx)_n=x_{n+1}$.
The system $(X,f)$ is called a {\em  mixing  subshift of finite type}. 
We fix $\nu \in (0,1)$ and consider the metric 
\vskip.15cm
$\hskip1.5cm \dist(x,y) = d_\nu(x,y)=\nu^{n(x,y)},
\;\text{ where }\;n(x,y)=\min\,\{ \,|i|\,: \; x_i \ne y_i  \}.
$
\vskip.15cm
\noindent The following sets play the role of the local stable and unstable 
manifolds of $x$:
$$
W^s_{\text{loc}}(x)=\{\,y\; | \;\, x_i=y_i, \;\;i\ge 0\,\}, \quad 
W^u_{\text{loc}}(x)=\{\,y\; | \;\, x_i=y_i, \;\;i\le 0\,\}
$$
Indeed, for all $x\in X$ and $n\in \N$,
$$
\begin{aligned}
&\dist (f^n x, f^n y )= \nu^n \cdot \dist  (x,y)  \quad\text{for all } y\in W^{s}_{\text{loc}}(x),\\
&\dist  (f^{-n}x, f^{-n}y )= \nu^n\cdot \dist  (x,y) \quad\text{for all } y\in W^{u}_{\text{loc}}(x),
\end{aligned}
$$
and for any $x, z\in X$ with $\dist(x,z) < 1$ the intersection of $W^s_{\text{loc}}(x)$ and $W^u_{\text{loc}}(z)$
consists of a single point. Thus, in this case we can take 
$\nu(x)=\nu$ and $\hat\nu (x) =\nu^{-1}$.

\vskip.2cm





\section{Proofs of Proposition \ref{bunching per} and Corollary \ref{bunching conj} } 

\subsection{Proof of Proposition \ref{bunching per}}
We consider the sequence of real-valued functions 
\begin{equation}\label{an}
  a_n(x)= \log \left( Q_\A(x,n)\cdot (\nu^n_x)^\beta \right)
  =\log \,\| \A_x^n \| + \log\, \| (\A_x^n)^{-1}\| +\log  (\nu^n_x)^\beta .
\end{equation}
It is easy to verify that this sequence of functions is {\em subadditive}, i.e.
$$
  a_{n+k}(x) \le a_k(x) +  a_n(f^k x) 
\quad\text{for all }x\in \M \text{ and  }n,k\in \N.
$$

We recall some results on subadditive sequences.
Let $f$ be a homeomorphism of a compact metric space $X$
and let $a_n$ be a subadditive sequence of continuous functions from $X$ to $\R$.
For  an ergodic $f$-invariant Borel  probability measure $\mu$  on $X$,  let
$
a_n(\mu)= \int_X a_n(x) \,d\mu.
$
The sequence 
of numbers $a_n (\mu) $ is subadditive, i.e.  $a_{n+k}(\mu) \le a_n(\mu) + a_k(\mu)$,\, and it is well known that 
$$
\lim_{n\to \infty} a_n (\mu) /n
 = \inf \{\, a_n (\mu)/ n\, :\; n\in \N \} =: \chi(a,\mu).
$$
By the Subaddititive Ergodic Theorem,  
$$
 \lim_{n\to \infty} a_n (x)/n= \chi(a,\mu) \quad
  \text{for $\mu$-almost all } x\in X.
$$ 

\begin{lemma}\cite[Proposition 4.9]{KS13} \label{subadditive}
Let $f$ be a homeomorphism of a compact metric space $X$
and  let $a_n : X \to \R$  be a subadditive sequence of continuous functions.

If $\chi(a,\mu)<0$  for every ergodic invariant Borel probability measure 
$\mu$ for $f$, then there exists $N$ such that $a_N(x) <0$ for all $x \in X$.
\end{lemma}

We will show that the assumption of the lemma is satisfied for the
sequence $a_n(x)$ given by \eqref{an}. We observe  that for this sequence,
$\,\chi(a,\mu)$
can be written in terms of Lyapunov exponents of the cocycle $\A$ and $\nu^\beta$.

\begin{definition}
Let $\mu$ be an ergodic $f$-invariant Borel probability measure on $\M$. The {\em upper and lower Lyapunov exponents} of $\A$ with respect to $\mu$ are
$$
\lambda_+(\A,\mu)= \lim_{n \to \infty} \frac 1n \log \| \A_x ^n \| \quad\text{and}\quad
\lambda_- (\A,\mu) = \lim_{n \to \infty} \frac 1n \log \| (\A_x ^n)^{-1} \|^{-1} 
$$
\end{definition}

\noindent By the Subadditive Ergodic Theorem, each of these limits exists 
and is the same $\mu$ almost everywhere.
 Now for the sequence in \eqref{an} for $\mu$ almost every $x$ we have:
\begin{equation}\label{chi}
\chi(a,\mu) = \lim_{n\to \infty} a_n (x)/n = \lambda_+(\A,\mu)- 
\lambda_- (\A,\mu) + \lambda (\nu^\beta, \mu),
\end{equation}
where $\lambda (\nu^\beta, \mu)=\lambda_+ (\nu^\beta, \mu)=\lambda_- (\nu^\beta, \mu)$ for the scalar cocycle $\nu^\beta$.
\vskip.1cm

The next theorem gives an approximation of the upper and lower Lyapunov 
exponents of a cocycle in terms of its periodic data. 
It is Theorem 1.4 and Remark 1.5 in \cite{KS16} stated for our setting.

\begin{theorem}\label{approx}
Let $(X,f)$ be a hyperbolic system, let $\mu$ be an ergodic $f$-invariant Borel probability measure on $X$, and let $\A$ be a H\"older continuous Banach 
cocycle over $f$.
Then for each  $\e>0$ there exists a periodic point $p=f^kp$ in $X$
such that 
\begin{equation}\label{la th}
\left| \,\lambda_+(\A,\mu)-  \frac1k \log \| \A_p ^k \| \, \right|<\e \quad\text{and}\quad
\left| \, \lambda_- (\A,\mu)- \frac1k \log \| (\A_p ^k)^{-1} \| ^{-1} \, \right|<\e.
\end{equation}

\noindent Moreover,  given finitely many  H\"older continuous Banach 
cocycle over $f$  and $K\in \N$, there exists a periodic point $p=f^kp$  with $k>K$
which gives simultaneous approximation \eqref{la th} for all the cocycles.
\end{theorem}

Let $\tilde L$ and $\tilde \theta <1$ be as in the assumption \eqref{fiber bunched per}.
We choose $ K\in \N$ such that $-3\e:=(\log \tilde L)/K + \log \tilde \theta <0$.  
Then for each  point
$p=f^kp$ with $k\ge K$ we have
$$
\frac1k \log \left(  \|\A_p^k\| \cdot \| (\A_p ^k)^{-1} \|\cdot (\nu_p^k)^\beta\right) \le \frac1k \log(\tilde L \, \tilde \theta^k)=
\,\frac1k \log \tilde L + \log \tilde \theta\,<-3\e.
$$
By Theorem \ref{approx}  there exists a periodic point $p=f^kp$ with  $k>K$ such that 
\eqref{la th} holds for $\A$, and also
for the scalar cocycle $\nu^\beta$ we have that 
$\,\la (\nu^\beta, \mu_p) = \frac1k \log (\nu_p^k)^\beta\,$ is $\e$ close to $\la (\nu^\beta, \mu)$.
Then by  \eqref{chi}  we have 
$$
\left| \,\chi(a,\mu) - \frac1k \log 
\left( \|\A_p^k\| \cdot \| (\A_p ^k)^{-1} \| \cdot (\nu_p^k)^\beta \right) \,\right|<3\e
$$
and hence  $\chi(a,\mu)<0$. Then by Lemma \ref{subadditive} there exists 
$N$ such that $a_N(x)<0$ for all $x$, and so
$\,Q_\A(x,N)\cdot (\nu^N_x)^\beta <1\,
\text{ for all }x\in \M.
$
By continuity, there exists $\theta<1$
such that the left hand side  is smaller than $\theta$ for
all $x$.
Writing  $n\in \N$  as $n=mN+r$, where $m\in \N\cup \{0\}$ and $0\le r<N$, 
we obtain
$$
 Q_\A(x,n)\cdot (\nu^n_x)^\beta \le  Q_\A(x,mN)\cdot (\nu^{mN}_x)^\beta
\cdot  Q_\A(f^{mN}x,r)\cdot (\nu^r_{f^{mN}x})^\beta \le L\, \theta^{\,n},
$$
 where  $L=\max \,\{\, Q_\A(x,r)\, (\nu^r_x)^\beta\,\theta^{-N}: \; x\in X, \;0\le r<N \}.$
The  inequality  with $\hat\nu$  
is obtained similarly, and 
we conclude that the cocycle $\A$ is fiber bunched.
$\QED$

 \vskip.2cm
 
\subsection{Proof of Corollary \ref{bunching conj}}

Since the cocycle $\A$ is fiber bunched, there exist numbers $L$ and $\theta<1$ 
such that $\,Q_\A(x,n) \cdot (\nu^n_x)^\beta < L\,\theta^n\,$
for all  $x\in \M$  and $n\in \N.$
It follows that whenever $p=f^kp$,
$$
  Q_\B(p,k) \cdot (\nu^k_p)^\beta \le \|C(p)\| \cdot Q_\A(p,k)\cdot  \|C(p)^{-1}\|
   \cdot (\nu^k_p)^\beta \le M^2 L\,\theta^k.
$$
Hence by Proposition \ref{bunching per} the cocycle $\B$ is fiber bunched.
$\QED$.


\section{Proof of Theorem \ref{conjugate data}}

\subsection{Cocycles over systems with a fixed point} 
The following result was  established  in \cite{S15} 
for $GL(d,\R)$-valued cocycles. 

\begin{theorem} \cite[Theorem 2.4]{S15} \label{old}
Suppose that $\A$ is fiber bunched and $\B$ has conjugate periodic data.
In addition, suppose that $f$ has a fixed point $p_0$  and the conjugacy $C(p)$ 
is $\beta$-H\"older continuous at  $p_0$,
i.e.  $d(C(p), C(p_0)) 
 \le c\, \dist (p,p_0)^\beta$ for every periodic point $p$.
Then there exists a unique $\beta$-H\"older continuous conjugacy 
$\bar  C$ between $\A$ and $\B$ such that $\bar  C(p_0)= C(p_0)$.
\end{theorem}

For Banach cocycles, the proof holds without modifications, except for using
Corollary  \ref{bunching conj}  instead of Proposition 2.3 in \cite{S15}
to obtain fiber bunching of $\B$.

We give an outline of the proof of this theorem since we will refer to it later. 
 We consider 
the cocycle $\tilde \B =C(p_0)\circ \B \circ C(p_0)^{-1}$ and  
the function $\tilde C(p)=C(p) C(p_0)^{-1}$,
so that $\tilde \B_{p_0} =\A_{p_0}$ and  $\tilde C(p_0)=\Id$.  
If $\tilde C(x)$ is a conjugacy  between $\A$ and $\tilde \B$ with $\tilde C(p_0)=\Id$,
then $ \bar C(x)= \tilde C(x) C(p_0)$ is a conjugacy between $\A$ and $\B$
with $\bar C(p_0)=C(p_0)$. 

Thus it suffices to consider the case when $\B_{p_0} =\A_{p_0}$
and $C(p_0)=\Id$.
First the conjugacy $\bar C(x)$ is constructed along 
the stable and unstable manifolds of $p_0$ using the stable and unstable holonomies 
(see Section \ref{holonomies sec} below):
\begin{equation}\label{Cs}
\begin{aligned}
 &  \bar C^s(x)=H^{\A,s}_{p_0,\,x}\circ H^{\B,s}_{x,\,p_0} \quad\text{for }x\in W^s(p_0),\\
 &\bar C^u(x)=H^{\A,u}_{p_0,\,x}\circ H^{\B,u}_{x,\,p_0} \quad\text{for }x\in W^u(p_0).
 \end{aligned}
\end{equation}
For a homoclinic point $x$ for $p_0$, i.e. $x\in W^s(p_0)\cap W^u(p_0)$,
it is shown that
$$\bar C^s(x)=\bar C^u(x)\overset{def}= \bar C(x).$$
 Thus we obtain $\bar C$ on the set of 
homoclinic points of $p_0$, which is dense in $\M$.
The function  $\bar C$ is $\beta$-H\"older continuous on this set, and hence
it can be extended to $\M$.
\vskip.2cm

\subsection{Holonomies} \label{holonomies sec} 
An important role in the proof of Theorem \ref{old}, as well as in the proof of 
Theorem \ref{conjugate data}, is played by holonomies. 
Holonomies were introduced in \cite{V08} and further studied in 
\cite{ASV,KS13,S15}. The result below gives existence of holonomies 
under the weakest fiber bunching assumption \eqref{fiber bunched}.

Let $\E=\M\times V$ be a trivial vector bundle over $\M$. 
For a Banach cocycle $\A$, 
we view $\A_x$ as a linear map from $\V_x$, the fiber at $x$, 
to $\V_{fx}$,\, so  $\,\A_x^n: \V_x \to \V_{f^nx}\,$ and 
$\,\A_x^{-n}: \V_x \to \V_{f^{-n}x}$.

\begin{proposition} \cite[Proposition 4.4]{S15} \label{existence of holonomies} 
Suppose that a cocycle $\A$ is fiber bunched. 
Then for every $x\in \M$ and $y\in W^s(x)$ the limit 
\begin{equation}\label{hol def}
 H^{\A,s}_{x,y} =\underset{n\to\infty}{\lim} \,(\A^n_y)^{-1} \circ \A^n_x, 
\end{equation}
exists and satisfies 
\begin{itemize}
\item[(H1)] $H^{\A,s}_{x,\,y}$ is a linear map from $\E_x$ to $\E_y$;
\vskip.1cm
\item[(H2)]  $H^{\A,s}_{x,\,x}=\Id\,$ and $\,H^{\A,s}_{y,\,z} \circ H^{\A,s}_{x,\,y}=H^{\A,s}_{x,\,z}$,\,\,
which implies $(H^{\A,s}_{x,\,y})^{-1}=H^{\A,s}_{y,\,x};$
\vskip.1cm
\item[(H3)]  $H^{\A,s}_{x,\,y}= (\A^n_y)^{-1}\circ H^{\A,s}_{f^nx ,\,f^ny} \circ \A^n_x\;$ 
for all $n\in \N$;
\vskip.1cm
\item[(H4)] $\| H^{\A,s}_{\,x,y} - \Id \,\| \leq c\,\dist (x,y)^{\beta},$
 where $c$ is independent of $x$   and $y\in W^{s}_{\text{loc}}(x).$

\end{itemize}
   \end{proposition}

The continuous map $H^{\A,s}:\;(x,y)\mapsto H^{\A,s}_{x,\, y}$,\, 
where $x\in \M$, $y\in W^s(x)$,  is called the (standard) {\em stable holonomy} 
for $\A$. The unstable holonomy $H^{\A,u}$ is defined 
similarly: 
$$H^{\A,u}_{x,\,y}\,=
\underset{n\to\infty}{\lim} \left( (\A^{-n}_y)^{-1} \circ (\A^{-n}_x) \right)
\,=
\underset{n\to\infty}{\lim} 
\left(\A^n_{f^{-n}y} \circ (\A^{n}_{f^{-n}x})^{-1} \right), \text{ where } y\in W^u(x).
$$
It satisfies (H1,\,2,\,4) and
\vskip.1cm

(H3$'$) $H^{\A,u}_{x,\,y}=
(\A^{-n}_y)^{-1} \circ H^{\A,u}_{f^{-n}x ,\,f^{-n}y}  \circ \A^{-n}_x\;$ 
 for all $n\in \N$.

\vskip.2cm

\subsection{Removing the fixed point assumption}
Now we obtain the H\"older continuous conjugacy  between the cocycles
assuming H\"older continuity of $C(p)$ at a periodic point $p_0$.
As was explained before,
we can assume that $\A_{p_0}=\B_{p_0}$ and $C(p_0)=\Id$.
Let $k$ be a period of the point $p_0$. Then $p_0$ is a fixed point for $f^k$
and so Theorem \ref{old} gives a unique H\"older conjugacy $\bar C$ between 
the iterated cocycles $\A^k$ and $\B^k$  over $f^k$ with $\bar C(p_0)=\Id$.
We will show that $\bar C$ is also a conjugacy between $\A$ and $\B$
following the approach of \cite{BK, Sch}.
\vskip.1cm

To simplify the notations we write $p$ for $p_0$ and $C$ for $\bar C$.
We consider a point $x\in W^s(p) \cap W^u(f^{k-1}p)$, so that $fx\in  W^u(p)$.
We will show that
\begin{equation}\label{Ax}
  \A_x=C(fx)\circ \B_x\circ C(x)^{-1}, \quad \text{i.e.\quad }
  C(x)= \A_x ^{-1} \circ C(fx)\circ \B_x.
\end{equation}
Since such points $x$ are dense in $\M$ and all the functions are continuous,
it follows that the equation holds for all $x\in \M$.

Now we prove \eqref{Ax}. Using equations \eqref{Cs} we obtain 
$$
\begin{aligned}
C(x) &= C^s(x) =H^{\A^k,s}_{p,\,x}\circ H^{\B^k,s}_{x,\,p} = 
\\ & =\underset{n\to\infty}{\lim} \,(\A^{nk}_x)^{-1} \circ \A^{nk}_p \circ  (\B^{nk}_p)^{-1} \circ \B^{nk}_x
=\underset{n\to\infty}{\lim} \, (\A^{nk}_x)^{-1} \circ  \B^{nk}_x, \quad\text{and}\\
 C(fx) &=  C^{u}(fx)=H^{\A^k,u}_{p,\,fx}\circ H^{\B^k,u}_{fx,\,p} =
\underset{m\to\infty}{\lim} \,\A^{mk}_{f^{-mk}fx} \circ (\B^{mk}_{f^{-mk}fx})^{-1}  \\
&= \underset{m\to\infty}{\lim} \,\A^{mk}_{f^{-mk+1}x} \circ (\B^{mk}_{f^{-mk+1}x})^{-1}.
 \end{aligned}
$$
Since $\A_x ^{-1} \circ C(fx)\circ \B_x$ equals 
$$
    \underset{m\to\infty}{\lim} 
   \,\A_x ^{-1} \circ \A^{mk}_{f^{-mk+1}x} \circ (\B^{mk}_{f^{-mk+1}x})^{-1}  \circ \B_x = 
   \underset{m\to\infty}{\lim} 
   \, \A^{mk-1}_{f^{-mk+1}x} \circ (\B^{mk-1}_{f^{-mk+1}x})^{-1},
    $$
we need to show  that 
$$
\underset{n\to\infty}{\lim} (\A^{nk}_x)^{-1} \circ  \B^{nk}_x  = \underset{m\to\infty}{\lim} 
   \, \A^{mk-1}_{f^{-mk+1}x} \circ (\B^{mk-1}_{f^{-mk+1}x})^{-1}. 
$$
As both limits exist, so does the following limit and it suffices to prove the equality
   $$
 \underset{m,n\to\infty}{\lim}  \, (\A^{nk}_x)^{-1} \circ  \B^{nk}_x \circ \B^{mk-1}_{f^{-mk+1}x} 
   \circ (\A^{mk-1}_{f^{-mk+1}x})^{-1} = \Id. 
$$

Since $x\in W^s(p)$ and $fx \in W^u(p)$, 
$$
f^{nk}x\to p \,\text{ as }n\to\infty  \quad\text{and}\quad  f^{-mk+1}x=f^{-mk}fx \to p 
\,\text{ as }m\to\infty.
$$
Moreover, by \eqref{Anosov def} 
there is a constant $c_1(x)$ such that for all $n,m\in \N$,
$$ 
\begin{aligned}
&\dist (f^{nk}x,\, p)< \nu^{nk}_x \cdot c_1(x)\, \dist_{W^s(p)}(x,p) = : c_2(x)\,\nu^{nk}_x \\
&\dist (f^{-mk+1}x, \,p)< \hat\nu^{-mk+1}_x\cdot c_1(x)\, \dist_{W^u(p)}(x,p) 
=:  c_3(x)\,\hat\nu^{-mk+1}_x.
\end{aligned}
$$
Let $\delta_0$ be as in Anosov Closing Lemma  \cite[Theorem 6.4.15]{KH}.
We take $\delta<\delta_0$ and let $m$ and $n$ be the smallest positive integers 
such both distances above are less than $\delta/2$. Then we have
\begin{equation}\label{delta}
\dist (f^{nk}x,\, f^{-mk+1}x) <\delta \quad\text{and}\quad
  \delta < c_4(x) \min \{\nu^{nk}_x,\hat\nu^{-mk+1}_x \} .
\end{equation}
Applying the  lemma to the orbit segment 
$
\,\{f^i x: \;\; i=-mk+1, \dots , nk\},\,
$
we obtain a periodic point $q=f^{(m+n)k-1}q\,$ such that 
$$
 \dist (f^i x,\, f^i q)\le L \delta \quad \text{for }i=-mk+1,\, \dots, nk.
$$
Let $z$ be the intersection point of  $W^s_{\text{loc}}(q)$ and  $W^u_{\text{loc}}(x)$. Then 
by the local product structure,
\begin{equation}\label{close}
\dist (f^i z,\, f^i x)\le c_5 \delta \,\text{ and }\,
\dist (f^i z,\, f^i q)\le c_5  \delta \quad \text{for }i=-mk+1, \dots , nk.
\end{equation}
Since  $z\in W^s(q)$, property (H3)  of the holonomies yields
$$
\A^{nk}_{z}=H^{\A,s}_{f^{nk}q,\,f^{nk}z}\circ \A^{nk}_{q} 
\circ H^{\A,s}_{z,\,q}, 
$$
and since $f^{nk} z\in W^u(f^{nk} x)$ using property  (H3$'$) we obtain
$$
\A^{nk}_{x}=  
H^{\A,u}_{f^{nk}z,\,f^{nk}x}\circ \A^{nk}_{z} \circ H^{\A,u}_{x,\,z}.
$$
Thus
$$
\A^{nk}_{x}=  
H^{\A,u}_{f^{nk}z,\,f^{nk}x}\circ H^{\A,s}_{f^{nk}q,\,f^{nk}z}\circ \A^{nk}_{q} 
\circ H^{\A,s}_{z,\,q} \circ H^{\A,u}_{x,\,z}.
$$
It follows from property (H4) that 
$$
H^{s, \A}_{z,q}=\Id+R^{s, \A}_{z,q}, \quad\text{where }\;
\|R^{s, \A}_{z,q}\| \le c \,\dist(z,q)^\beta \le c_6  \delta^\beta,
$$
Similar estimates hold for the other holonomies, as well as their inverses,  due to \eqref{close}.
Thus we obtain 
\begin{equation}\label{A}
(\A^{nk}_{x})^{-1}=(\Id +R_1)\circ (\A^{nk}_{q})^{-1} \circ (\Id +R_2), \quad\text{where }\;
\|R_1\|, \|R_2\|  \le c_6 \delta^\beta,
\end{equation}
and similarly, 
\begin{equation}\label{A2}
(\A^{mk-1}_{f^{-mk+1}x})^{-1} =
(\Id +R_3)\circ (\A^{mk-1}_{f^{-mk+1}q})^{-1}  \circ (\Id +R_4), 
\end{equation}
where $\|R_3\|, \|R_4\|  \le c_6 \delta^\beta$. \,Let  $q'=f^{-mk+1}q=f^{nk}q$ and $x'=f^{-mk+1}x$.
Then 
\begin{equation}\label{Bx}
\B^{nk}_x \circ \B^{mk-1}_{f^{-mk+1}x} = \B^{(m+n)k-1}_{x'}=
(\Id +R_5)\circ \B^{(m+n)k-1}_{q'} \circ (\Id +R_6), 
\end{equation}
 where 
$\|R_5\|, \|R_6\|  \le c_6 \delta^\beta.$
 Since $q'$ is a periodic point of period $(m+n)k-1,\,$ by the assumption 
there exists $C(q')$ such that 
\begin{equation} \label{q'}
\begin{aligned}
 & \B^{(m+n)k-1}_{q'}=C(q') \circ \A^{(m+n)k-1}_{q'} \circ C(q')^{-1},
 \;\text{ where } C(q')=\Id+R_7, \\
 & C(q')^{-1}=\Id+R_8 \quad\text{with }
 \|R_7^n\|, \; \|R_8^n\| \le c_7'\, \dist (p,q')^\beta  \le c_7 \delta^\beta.
 \end{aligned}
\end{equation}
It follows from  \eqref{Bx} and \eqref{q'} that
\begin{equation}\label{B2}
\B^{nk}_x \circ \B^{mk-1}_{f^{-mk+1}x} =(\Id + R_9)\circ  \A^{nk}_{q} \circ \A^{mk-1}_{f^{-mk+1}q} 
\circ (\Id +R_{10}), 
\end{equation}
where $\|R_9\|, \|R_{10}\|  \le c_8 \delta^\beta.$
\vskip.1cm

 Using \eqref{A}, \eqref{A2}, and \eqref{B2}
 and combining terms of type $\Id+R_i$ we obtain 
$$
\begin{aligned}
& (\A^{nk}_x)^{-1} \circ\, \B^{nk}_x \circ \B^{mk-1}_{f^{-mk+1}x} \circ  (\A^{mk-1}_{f^{-mk+1}x}) ^{-1}  =\\
& (\Id +R_1)  \circ (\A^{nk}_{q})^{-1} \circ (\Id +R_{11}) 
 \circ \A^{nk}_{q} \circ \A^{mk-1}_{f^{-mk+1}q} \circ  
 (\Id +R_{12}) \circ (\A^{mk-1}_{f^{-mk+1}q})^{-1} \circ (\Id +R_{4}) \\
& = \Id +  (\A^{nk}_{q})^{-1} \circ R_{11} \circ  \A^{nk}_{q}  \circ \A^{mk-1}_{f^{-mk+1}q} 
\circ  R_{12} \circ (\A^{mk-1}_{f^{-mk+1}q})^{-1}
\end{aligned}
$$
$$
+ (\A^{nk}_{q})^{-1} \circ R_{11}\circ  \A^{nk}_{q}  + \A^{mk-1}_{f^{-mk+1}q} 
\circ  R_{12} \circ (\A^{mk-1}_{f^{-mk+1}q})^{-1}+ \,\,\text{smaller terms}.\;\;
$$
Since $\|R_i \|\le  c_9 \delta^\beta$ for each $i$, and $\delta$ satisfies \eqref{delta},
and the cocycle is fiber bunched, we have
$$
\begin{aligned}
 & \|\A^{nk}_{q}\| \cdot  \|(\A^{nk}_{q})^{-1} \| \cdot \|R_{11}\| \le 
 Q_\A (q,nk)\cdot c\delta^\beta \le  Q_\A (q,nk)\cdot c_{10}(\nu_x^{nk})^\beta =
 c_{10} L\theta^{nk} \quad\text{and}
 \\
& \| \A^{mk-1}_{f^{-mk+1}q} \|\cdot \| (\A^{mk-1}_{f^{-mk+1}q})^{-1}\|
\cdot \| R_{12} \|  \le 
c_{10}(\hat \nu_x^{-mk+1})^\beta = c_{10} L\theta^{mk-1}.
\end{aligned}
$$
Thus 
   $$
   \lim_{n,m\to\infty} (\A^{nk}_x)^{-1} \circ  \B^{nk}_x \circ \B^{mk-1}_{f^{-mk+1}x} 
   \circ (\A^{mk-1}_{f^{-mk+1}x})^{-1} =\Id,
   $$
   and so \eqref{Ax} holds.  
   $\QED$


\section{Proof of Theorem \ref{measurable conjugacy}}
We consider the strong operator topology, i.e. the topology of pointwise convergence, 
on the space of linear operators $L(V)$. It is  induced by the family of semi-norms 
$P_v(A)=\|A(v)\|$, where $v\in V$. 

Since $V$ is separable, this  topology is separable and  metrizable
on any  set $G\subset L(V)$ that is bounded in norm. Indeed, let $\{v_n: n\in \N\}$ 
be a countable dense set in the unit ball in $V$. Then 
\begin{equation}\label{bar d}
  \bar d (A,B)=\,\sum_{n=1}^\infty 
  \,\frac{\|A(v_n)-B(v_n)\|}{1+\|A(v_n)-B(v_n)\|}\cdot 2^{-n}
\end{equation}
is a distance on $L(V)$. The convergence in $\bar d$ is the pointwise convergence on
the set $\{v_n: n\in \N\}$.  It induces the strong operator topology on $G$ since for a bounded 
sequence convergence on each $v_n$ is equivalent to convergence on each $v\in V$.
To  show separability, we take  the set $\{v_n\}$ in the definition of $\bar d$ to be linearly 
independent and consider a  countable dense set $U=\{u_n: n\in \N\}$ in $V$. Then the 
set of all finite sequences $\{u_{n_1}, \dots ,u_{n_k} \}$ in $U$ is countable and for each 
finite sequence we can take a bounded operator $B_{n_1,\dots , n_k}$ 
such that 
$$
B_{n_1,\dots , n_k}(v_i)=u_{n_i} \,\text{ for  }i=1, \dots , k. 
$$
This can be done by extending the coordinate functionals $x_i$ on 
span$\{ v_1, \dots , v_k\}$ to $V$ using Hahn-Banach Theorem and defining
$
B_{n_1,\dots , n_k}(v)=\sum_{i=1}^k x_i(v) u_{n_i}.
$
The set of such operators is a countable dense set in $(L(V), \bar d)$.
Indeed, given $\e>0$ and $A\in L(V)$ we can fix a large $k$ so that the ``tail"
of the series in \eqref{bar d} is small
and then choose $u_{n_i}$ sufficiently close to $A(v_i)$
 so that 
 $$
 \|A(v_i)-B_{n_1,\dots , n_k}(v_i)\|=\|A(v_i)-u_{n_i}\| \quad \text{is small for  }\,i=1, \dots ,k.
 $$
Thus the strong operator  topology is separable and hence second countable on $G$.
 \vskip.05cm

The corresponding strong operator topology on $GL(V)$ is induced by the 
embedding $i: GL(V) \to L(V) \times L(V)$ given by $i(A)=(A,A^{-1})$.
The convergence in this topology is the pointwise convergence
of operators and their inverses: a sequence $A_n$ converges to $A\,$ 
 if for each $v\in V$,
$$
  \|A_n(v)-A(v)\|+\|A_n^{-1}(v)-A^{-1}(v)\| \to 0 \quad\text{as }n\to\infty.
$$ 
It follows from the results for $L(V)$ that  this topology is separable and  metrizable by
$$
  \tilde d (A,B)= \bar d (A,B) +  \bar d (A^{-1},B^{-1})
$$
on any set $G\subset GL(V)$ bounded with respect to the metric $d$ given by \eqref{d}.

By the Uniform Boundedness Principle, a sequence $A_n$ 
 that converges in strong operator topology on $GL(V)$ is bounded in the metric $d$,
 and hence any subset of $GL(V)$ that is compact in the strong operator topology
 is bounded in $d$.
Since the cocycle $\B$ takes values in such a subset, $\B$  is uniformly bounded, 
and thus satisfies the fiber bunching condition \eqref{fiber bunched}. 
\vskip.1cm

Let $C$ be a $\mu$-measurable conjugacy between $\A$ and $\B$.
First we show that $C$ intertwines holonomies of $\A$ and $\B$
on a set of full measure, i.e.
there exists a set $Y\subset \M$ with $\,\mu(Y)=1$  such that 
\begin{equation}\label{int meas}
H_{x,y}^{\A,s}=C(y)\circ H_{x,y}^{\B,s} \circ C(x)^{-1}
\quad\text{for all }x,y\in Y\; \text{ such that }y\in W^s(x),
\end{equation}
and a similar statement holds for the unstable holonomies. Since 
$$
C(x)= (\A^n_x)^{-1}\circ C(f^n x)\circ \B^n_x \quad\text{and}\quad
 H^{\A,s}_{x,\,y}= (\A^n_y)^{-1}\circ H^{\A,s}_{f^nx ,\,f^ny} \circ \A^n_x,
$$
it suffices to prove that 
$\,H_{f^nx,f^ny}^{\A,s}=C(f^ny)\circ H_{f^nx,f^ny}^{\B,s} \circ C(f^nx)^{-1}.\,$
Thus we can assume that $y$ lies on the local stable manifold of $x$.
\vskip.1cm

Let $x\in \M$ and  $y\in W^s_{\text{loc}} (x)$.
Since $\A_x=C(fx) \circ \B_x \circ C(x)^{-1}$,
we have
\begin{equation}\label{AA}
\begin{aligned}
&(\A^n_y)^{-1} \circ  \A^n_x 
\,=\, C(y) \circ (\B^n_y)^{-1} \circ C(f^n y)^{-1} \circ  C(f^nx) \circ \B^n_x
\circ  C(x)^{-1}= \\
&= C(y)\circ (\B^n_y)^{-1} \circ (\Id +r_n) \circ \B^n_x\circ C(x)^{-1}= \\
&= C(y) \circ (\B^n_y)^{-1} \circ \B^n_x\circ C(x)^{-1} + 
C(y)\circ (\B^n_y)^{-1} \circ r_n \circ \B^n_x \circ C(x)^{-1}, 
\end{aligned}
\end{equation}
where 
$$
 r_n  = C(f^n y)^{-1} \circ  C(f^nx)-\Id
=  C(f^n y)^{-1} \circ ( C(f^nx)- C(f^ny)).
$$

Since $C$ is $\mu$-measurable, $\|C\|$ and $\|C^{-1}\|$ are measurable 
functions from $\M$ to $\R$ and hence  there exists a 
compact set $S_1\subset \M$ with $\mu(S_1)>3/4$ such that $\|C\|$
and $\|C^{-1}\|$ are bounded on $S_1$.
Let  $G\subset GL(V)$ be  a $d$-bounded set that contains 
the values of $C$ and $C^{-1}$ on the set $S_1$.
Then $C$ is a  $\mu$-measurable function from $S_1$ to 
the separable, and hence second countable,  metric space $(G,\tilde d)$.
Hence by Lusin's theorem there exists a 
compact set $S\subset S_1$ with $\mu(S)>1/2$ such that $C$
is uniformly continuous on $S$.

Let $Y$ be the set of points in $\M$ for which the frequency of 
visiting $S$ equals $\mu(S)>1/2$. By Birkhoff Ergodic Theorem,
$\mu(Y)=1$. If $x$ and $y$ are in $Y$, there exists a sequence $\{n_i\}$
such that $f^{n_i}x$ and $f^{n_i}y$ are in $S$ for all $i$.
Thus  $\|C\|$, $\|C^{-1}\|$ are uniformly bounded on the set $\{f^{n_i}x, \, f^{n_i}y\}$
and $\| (C(f^{n_i}x)- C(f^{n_i}y))(v)\|\to 0$ for every  $v\in V$.
Hence  
\begin{equation}\label{ni}
\text{the sequence $\|r_{n_i}\|$ is bounded\, and }\;
\| r_{n_i}(v) \|\to 0 \quad\text{for every  }v\in V.
\end{equation}
We rewrite the formula \eqref{AA} as 
\begin{equation}\label{AA2}
C(y)^{-1}  \circ (\A^n_y)^{-1} \circ  \A^n_x  \circ C(x) -
(\B^n_y)^{-1} \circ \B^n_x = (\B^n_y)^{-1} \circ r_n \circ \B^n_x. 
\end{equation}
By the definition \eqref{hol def}, the left hand side converges in the operator norm to 
\begin{equation}\label{HC}
C(y)^{-1}  \circ H^{\A,s}_{x,y} \circ C(x) -H^{\B,s}_{x,y}.
\end{equation}
Hence the right hand side of \eqref{AA2} has a limit in the operator norm.
To prove that it equals 0, it suffices to show  that for each $v\in V$,  
$\,((\B^n_y)^{-1} \circ r_n \circ \B^n_x)(v)$ tends to 0
along a subsequence.

Let $\{n_i\}$ be the sequence  in  \eqref{ni}.
Since $\B$ takes values in a compact subset of $GL(V)$ with strong operator topology,
there exists a subsequence $\{ n_{i_j}\}$ of $\{n_i\}$ such that 
$(\B^{n_{i_j}}_x)(v)$  converges 
for every $v$ and, as we observed above,
 $(\B_y^n)^{-1}$ is uniformly bounded. Now it follows  from \eqref{ni}
that
$$
\left( (\B^{n_{i_j}}_y)^{-1} \circ r_{n_{i_j}} \circ \B^{n_{i_j}}_x \right)(v) \to 0
  \quad\text{for every  }v\in V.
$$
We conclude that \eqref{HC} equals 0
and so \eqref{int meas} follows.
The statement for the unstable holonomies is proven similarly.
\vskip.1cm

Now we establish H\"older continuity of $C$ on a set of full measure.
We consider a small open set $U$ in $\M$ with the product structure of stable and unstable manifolds, 
$$
U=W^s_{\text{loc}}( x_0)\times W^u_{\text{loc}}( x_0) =\,
\{W^s_{\text{loc}}( x)\cap W^u_{\text{loc}}( y)\,: \; x\in W^s_{\text{loc}}( x_0) \text{ and } y\in W^u_{\text{loc}}( x_0)\}.
$$
We take a finite cover of $\M$ by such sets. It suffices to show that $C$  is H\"older continuous
on a full measure subset of each such set $U$.

Since the measure $\mu$ has local product structure,
$\mu$ is equivalent to the product of conditional measures on $W^s_{\text{loc}}( x_0)$
and $W^u_{\text{loc}}( x_0)$, and hence for $\mu$ almost every local 
stable leaf in $U$, the set of  points of 
$Y$ on the leaf has full conditional measure.  
Let $Y_U$ be the the set of points in $Y\cap U$ that lie on such leaves.
Then $Y_U$ has full measure in $U$.
Since the holonomies of the unstable foliation are absolutely continuous
with respect to the conditional measures, 
 for any two  points $x$ and $z$ in $Y_U$, there exists a point 
$y\in W^s_{\text{loc}}(x)\cap Y_U$  such that 
$y'=W^u_{\text{loc}}(y)\cap W^s_{\text{loc}}(z)$ is also in $Y_U$.
\vskip.1cm

Now we show that $C$ and $C^{-1}$ are bounded on $Y_U$.  We fix $x\in Y_U$ and for any 
$z\in Y_U$ we consider $y$ and $y'$ as above. Then by equation \eqref{int meas} and property (H4) 
of the holonomies  we have
\begin{equation}\label{Cxy}
C(y)= H_{x,y}^{\A,s} \circ C(x) \circ (H_{x,y}^{\B,s}) ^{-1} =  
(\Id +R^{\A,s}_{x,y}) \circ C(x)\circ (\Id +R^{\B,s}_{x,y}),
\end{equation}
where $\,\|R^{\A,s}_{x,y}\|, \; \|R^{\B,s}_{x,y}\| \le c\,\dist (x,y)^\beta.$
If the set $U$ is sufficiently small, then the terms $R$ are less than one in norm 
 and it follows that $\|C(y)\|\le 4 \|C(x)\|$.
Also considering the pairs $y,y'$ and $y',z$ we conclude that
$\|C(z)\|\le 4^3 \|C(x)\|$, and similarly  $\|C(z)^{-1}\| \le 4^3 \|C(x)^{-1}\|$.
\vskip.05cm

Now we establish  H\"older continuity of $C$ on $Y_U$. 
We fix $x,z\in Y_U$ and consider $y$ and $y'$ as above. 
Then by \eqref{Cxy} we have
$$
\begin{aligned}
& C(y) \circ C(x)^{-1} = (\Id +R^{\A,s}_{x,y} )\circ  C(x) \circ (\Id +R^{\B,s}_{x,y}) \circ C(x)^{-1}= \\
& =  \; \Id+R^{\A,s}_{x,y}+ C(x) \circ R^{\B,s}_{x,y} \circ C(x)^{-1} +
R^{\A,s}_{x,y} \circ  C(x) \circ R^{\B,s}_{x,y} \circ C(x)^{-1}.
\end{aligned}
$$
Since $C$ and $C^{-1}$  are bounded on $Y_U$, it follows that 
$$
 \|C(y) \circ C(x)^{-1} -\Id\|\le c_1\,\dist(x,y)^\beta.
$$
Hence we have
$$
\begin{aligned}
 & d(C(x),C(y))=\|C(x)-C(y)\|+ \|C(x)^{-1}-C(y)^{-1}\|\le\\
 &\le  \|C(x)C(y)^{-1}-\Id \|\cdot \|C(y)\|+
 \|C(x)^{-1}\|\cdot \|\Id-C(x)C(y)^{-1}\| \le \\ &\le c_2\, \dist(x,y)^\beta,
 \quad\text{where $c_2$ does not depend on $x$ and $y$.}
\end{aligned}
$$
Using similar estimates for $y,y'$ and $y',z$
and the local product structure of the stable and unstable manifolds 
we conclude that for all $x,z\in Y_U$,
$$
d(C(x), C(z))\le c_3 \,\dist(x,z)^\beta.
$$

Thus we obtain H\"older continuity of $C$ on a set of full measure
$Y_1 \subseteq Y$.
Let $Y_2 = \bigcap_{n=-\infty}^{\infty} f^n(Y_1)$. 
Then $Y_2$ is $f$-invariant and 
$\A(x)=C(fx)\circ \B(x) \circ C(x)^{-1}$ for all $x\in Y_2$. Since $\mu$ has
full support and $\mu(Y_2)=1$, the set $Y_2$ is dense in $\M$. 
Hence $C$ extends
from $Y_2$ to a H\"older continuous conjugacy
 $\tilde C$ on $\M$. 
 $\QED$




\end{document}